\documentclass[12pt, reqno]{amsart}

\usepackage{latexsym}
\usepackage{amssymb}
\usepackage{mathrsfs}
\usepackage{amsmath}
\usepackage{fancybox,color}
\usepackage{enumerate}
\usepackage[latin1]{inputenc}
\usepackage[active]{srcltx}
\usepackage[colorlinks]{hyperref} 
\usepackage{hypernat}
\usepackage{eurosym}
\usepackage{url}

\newcommand{\kom}[1]{}
\renewcommand{\kom}[1]{{\bf [#1]}}

\definecolor{gruen}{cmyk}{1.0,0.2,0.7,0.07}

 \def\1{\raisebox{2pt}{\rm{$\chi$}}}

\newtheorem{theorem}{Theorem}

\newtheorem{proposition}[theorem]{Proposition}

\newcommand{\R}{{\mathbb R}}

 \def\1{\raisebox{2pt}{\rm{$\chi$}}}
 
\newcommand{\abs}[1]{\left|#1\right|}

\newcommand{\Rn}{\mathbb{R}^n}

\def\vint_#1{\mathchoice%
          {\mathop{\kern 0.2em\vrule width 0.6em height 0.69678ex depth -0.58065ex
                  \kern -0.8em \intop}\nolimits_{\kern -0.4em#1}}%
          {\mathop{\kern 0.1em\vrule width 0.5em height 0.69678ex depth -0.60387ex
                  \kern -0.6em \intop}\nolimits_{#1}}%
          {\mathop{\kern 0.1em\vrule width 0.5em height 0.69678ex
              depth -0.60387ex
                  \kern -0.6em \intop}\nolimits_{#1}}%
          {\mathop{\kern 0.1em\vrule width 0.5em height 0.69678ex depth -0.60387ex
                  \kern -0.6em \intop}\nolimits_{#1}}}
\def\vintslides_#1{\mathchoice%
          {\mathop{\kern 0.1em\vrule width 0.5em height 0.697ex depth -0.581ex
                  \kern -0.6em \intop}\nolimits_{\kern -0.4em#1}}%
          {\mathop{\kern 0.1em\vrule width 0.3em height 0.697ex depth -0.604ex
                  \kern -0.4em \intop}\nolimits_{#1}}%
          {\mathop{\kern 0.1em\vrule width 0.3em height 0.697ex depth -0.604ex
                  \kern -0.4em \intop}\nolimits_{#1}}%
          {\mathop{\kern 0.1em\vrule width 0.3em height 0.697ex depth -0.604ex
                  \kern -0.4em \intop}\nolimits_{#1}}}

\newcommand{\aveint}[2]{\mathchoice%
          {\mathop{\kern 0.2em\vrule width 0.6em height 0.69678ex depth -0.58065ex
                  \kern -0.8em \intop}\nolimits_{\kern -0.45em#1}^{#2}}%
          {\mathop{\kern 0.1em\vrule width 0.5em height 0.69678ex depth -0.60387ex
                  \kern -0.6em \intop}\nolimits_{#1}^{#2}}%
          {\mathop{\kern 0.1em\vrule width 0.5em height 0.69678ex depth -0.60387ex
                  \kern -0.6em \intop}\nolimits_{#1}^{#2}}%
          {\mathop{\kern 0.1em\vrule width 0.5em height 0.69678ex depth -0.60387ex
                  \kern -0.6em \intop}\nolimits_{#1}^{#2}}}
\newcommand{\ud}{\, d}

\newcommand{\parts}[2]{\frac{\partial {#1}}{\partial {#2}}}
\newcommand{\ol}{\overline}
\newcommand{\Om}{\Omega}

\begin{document}

\title[Infinite initial values]{A remark on infinite initial values for quasilinear parabolic equations}

\author[Lindqvist]{Peter Lindqvist}
\address{Department of Mathematics and Statistics, Norwegian University of Science and Technology, NO-7491 Trondheim, Norway}
\email{lqvist@math.ntnu.no}

\author[Parviainen]{Mikko Parviainen}
\address{Department of Mathematics and Statistics, University of
Jyv\"askyl\"a, PO~Box~35, FI-40014 Jyv\"askyl\"a, Finland}
\email{mikko.j.parviainen@jyu.fi}

\date{\today}

\begin{abstract}
  We study the possibility of prescribing infinite initial values for solutions of the Evolutionary $p$-Laplace Equation in the fast diffusion case $p > 2$. This expository note has been extracted from our previous work. When infinite  values are prescribed on the whole initial surface,
   such solutions can exist only if the domain is a space-time cylinder. 
\end{abstract}

\maketitle

We are interested in a violent initial value problem for the Evolutionary $p$-Laplace Equation 
\begin{align*}
\parts{u}{t}=\operatorname{div}(\abs{\nabla u}^{p-2}\nabla u),\quad p>2,
\end{align*}
in a cylindrical domain $\Om_T=\Om\times(0,T)$, where $\Om$ is an open connected domain in the Euclidean $n$-dimensional space $\Rn$. We study the  slow diffusion case $p>2$. For $p=2$, the equation is the heat equation 
\begin{align*}
\parts{u}{t}=\Delta  u,
\end{align*}
but the phenomenon we describe is impossible for the heat equation.

The objective of our note is the possibility of prescribing  \emph{infinite} initial values in some subset $E$ of $\Om$, i.e.\
\begin{align}
\label{eq:infty-init}
\lim_{t\to 0+}u(x,t)=\infty \quad \text{ at every } x\in E.
\end{align}
We assume that $u\ge 0$ and do not mind the boundary values  at the lateral boundary $\partial \Om \times [0,T]$, since they now play a minor role in comparison to the violent effects of the $\infty $ initial values. The solutions are defined as weak solutions belonging to the Sobolev space $L^{p}_{\text{loc}}(0,T;W^{1,p}_{loc}(\Om))$ or even to $L^{p}_{\text{loc}}(0,T;W^{1,p}(\Om))$, and they are continuous in $\Om_T$. See for example  \cite{dibenedetto93} or \cite[Definition 5]{kuusilp16}. 

We have detected a strange phenomenon, which can be described as follows.
\begin{itemize}
\item If the $n$-dimensional Lebesgue measure  $\abs{E}>0$  of $E$ is strictly positive, then $E=\Om$.
\item There exist solutions, if $E=\Om$ and $\Om$ is bounded.
\item  There is no solution, if $\Om=\Rn$ and  $\abs{E}>0$.
\end{itemize}
To spell it out, the condition that
\begin{align*}
\liminf_{t\to 0+} u(x_0,t)<\infty
\end{align*}
at a point $x_0$ in $\Om$ is incompatible with the condition that (\ref{eq:infty-init}) holds in a set $E$ of positive $n$-dimensional measure. In the case when $\Om$ is the whole space $\Rn$, the \emph{global} restriction that $u\geq 0$ is decisive.\footnote{Otherwise one can construct  an odd periodic extension of $U$, first defined in the box $\Om \,=\,[0,1]\times\cdots\times [0,1]$, to the whole $\Rn$. The corresponding $V$ will do as an example.}

We remark that if $\Om$ is bounded and regular, then there exists a separable solution of the form 
\begin{align}
\label{eq:separable}
V(x,t)=t^{-\frac{1}{p-2}}U(x)
\end{align}
 where $U\in C(\ol \Om)$, $U>0$ in $\Om$ and $U|_{\partial \Om}=0$. The function $U$ solves the auxiliary equation  
 \begin{align*}
\operatorname{div}(|\nabla {U}|^{p-2}\nabla {U}\bigr)\:+\;\tfrac{1}{p-2}\, {U}\:=\;0.
\end{align*}
The function $V$ can serve as a minorant: every other solution $u\ge 0$ with $E=\Om$, the lateral boundary values of which are not necessarily zero, has to grow faster than $V$. That is 
 \begin{align*}
u(x,t)\ge t^{-\frac{1}{p-2}}U(x) \quad \text{in }\Om_T.
\end{align*}  
\begin{proposition}
If a weak solution $u\ge 0$, satisfies 
\begin{align*}
\lim_{t\to 0+} u(x,t)=\infty
\end{align*} 
at every point $x$ in  a set of positive $n$-dimensional measure, then the limit is infinite at each point $x\in \Om$. Moreover, 
\begin{align}
\label{eq:rate-of-infinity}
\liminf_{(x,t)\to (x_0,0)} \bigg[t^{\frac{1}{p-2}} u(x,t)\bigg]>0
\end{align}
at every point $x_0$ in $\Om$, and of course $(x,t)\in \Om_T$. 
\end{proposition}

The case $\abs{E}=0$ remains. For example, the celebrated Barenblatt solution for $t>0$ 
\begin{equation}
\label{eq:barenblatt}
\begin{split}
  {\mathcal B}(x,t)=t^{-n/\lambda}\Bigg\{C-\frac{p-2}{p}\lambda^{\frac{1}{1-p}}\Big(\frac{\abs x}{t^{1/\lambda}}\Big)^{\frac{p}{p-1}}\Bigg\}_+^{\frac{p-1}{p-2}}\,,
\end{split}
\end{equation}
where $\lambda=n(p-2)+ p$,  the constant $C>0$ can be chosen freely, and $\{\cdot\}_+$ means $\max\{\cdot,0\}$, exhibits that it is possible that $E$ consists of a single point: $E=\{  (0,0) \}$. We have that ${\mathcal B}(x,0)=c\delta_0$ (Dirac's delta) in the sense that 
\begin{align*}
\lim_{t\to 0+} \int_{\Rn} {\mathcal B}(x,t) \phi(x)\ud x=c\phi(0)
\end{align*}
for all test functions $\phi\in C^\infty_0(\Rn)$. In general, the question about whether there exists a Radon measure $\mu$ such that '$u=\mu$' as $t=0$ has drawn much attention, see for example \cite{aronsonc83, benilancp84, dahlbergk84} in the case of the porous medium equation.
It is clear that no $\sigma$-finite measure will do in the case $E=\Om$, since (\ref{eq:rate-of-infinity}) holds. However, if $\abs{E}=0$ this is possible.  The initial values can even be prescribed as a Radon measure, see  \cite{dibenedettoh89} and  \cite{kuusip09}.

Finally, we mention a surprising result, when the domain is of the form 
\begin{align*}
 \Om^{\Psi}\,=\,\{ (x,t) \,:\,x\in \Om,\ \Psi(x)<t<T\}
\end{align*}
 where $\Psi$ is a given smooth function. Assume now that
 \begin{align*}
\lim_{(y,\tau)\to (x,\Psi(x))} u(x,t)=+\infty
\end{align*}
at every point $x\in \Om$, where of course $(y,\tau)$ is in  $ \Om^{\Psi}$.  Indeed, such a solution does not exist, except when $\Psi(x)\equiv \text{constant}$. Then we are back to the case $\Om\times (0,T)$. \emph{Infinite  initial values cannot be prescribed on a whole surface in the $(x,t)$-space which is not a timeslice $t= constant$.}
 It is arresting, that not even a domain restricted by a hyperplane
\begin{align*}
t_0+a_1x_1+\ldots+a_nx_n<t<T,\,\, \abs{x}<1, 
\end{align*}
is possible, except when $a_1=\ldots=a_n=0$.

For the proof we let $x_0 \in \Om$
be a point such that a small ball $B=B(x_0,r)$ is divided by the smooth level surface  $C_0 =\{x\in \Om\,:\,\Psi(x) =  \Psi(x_0)\}$ into two pieces, say $B^+ = \{x\in  B\,:\,\Psi(x) > \Psi(x_0)\}$ and $B^- =  \{x\in B\,:\,\Psi(x) < \Psi(x_0)\}$. Then let $\omega = \omega(x)$ be the solution of the elliptic problem
\begin{equation*}
  \Delta_p\,\omega\,=\,0\quad \text{in}\quad B^-\quad \text{and}\quad   \omega =  \begin{cases}   0\quad\text{in}\quad B^-\cap\partial B\\  1\quad\text{in}\quad C_0\cap  B .\end{cases}
\end{equation*}
(If desired, we may round off the boundary values so that the value 1 is prescribed in a smaller set, yet one of positive area). On the parabolic boundary of  the domain
$$ \Psi(x) < t < \Psi(x_0),\quad x\in B^-$$
we have $u(x,t) \geq k\omega(x)$ for any constant $k$. By the comparison principle, see \cite {kilpelainenl96} and \cite[Theorem 2.4]{bjornbgp15}, the inequality
$u(x,t) \geq k\omega(x)$ holds in this domain. Sending $k$ to $\infty$ we obtain that $u \equiv \infty$ in an open set of $\R^{n+1}$. This contradicts the continuity of $u$.

  For the Porous medium equation 
 \begin{align*}
\parts{u}{t}=\Delta (u^m),\quad m>1,
\end{align*} 
 a similar theory is valid. The separable solution corresponding to (\ref{eq:separable}) is called 'the friendly giant' by Dahlberg and Kenig \cite{dahlbergk88}, see also Vazquez \cite[p.\ 111]{vazquez07}. 
 
As mentioned above, most of this note has been extracted from our work \cite{kuusilp16}. For more recent related works, see \cite{kinnunenl16, kinnunenllp, parviainenv}.

\def\cprime{$'$} \def\cprime{$'$} \def\cprime{$'$}

\end{document}